\documentclass[11pt]{article}
\newcommand{\proof}{\noindent {\bf Proof: }}
\newcommand{\remark}{\noindent {\bf Remark: }}

\newtheorem{theorem}{Theorem}
\newtheorem{lemma}{Lemma}

\newtheorem{proposition}{Proposition}
\newtheorem{defi}{Definition}
\def\qed{\hfill $\Box$}

\usepackage{amssymb}
\usepackage{graphicx}

\begin{document}

\title{Bounded representation and radial projections of bisectors in normed spaces}
\author{\'A. G. Horv\'ath and H. Martini}
\date{~}

\maketitle

\begin{abstract}
It is well known that the description of topological and geometric properties of bisectors in normed spaces is a non-trivial subject.
In this paper we introduce the concept of bounded representation of bisectors in finite dimensional real Banach spaces. This useful notion combines the concepts of bisector and shadow boundary of the unit ball, both corresponding with the same spatial direction. The bounded representation visualizes the connection between the topology of bisectors and shadow boundaries (Lemma 1) and gives the possibility to simplify and to extend some known results on radial projections of bisectors. Our main result (Theorem 1) says that in the manifold case the topology of the closed bisector and the topology of its bounded representation are the same; they are closed, $(n-1)$-dimensional balls embedded in Euclidean $n$-space in the standard way.
\end{abstract}

{\bf MSC(2000):} 46B20, 51M05, 52A21, 57N16

\textbf{Keywords: } bisector, bounded representation, Minkowski space, normed space, radial projection, shadow boundary, topological manifold

\section{Introduction and some preliminary results}

In recent times, Minkowski Geometry (i.e., the geometry of finite dimensional, real Banach spaces; see \cite{tho}) became again an important research field. Strongly related to Banach Space Theory, Finsler Geometry and Classical Convexity, it is permanently enriched by new results in the spirit of Discrete and Computational Geometry, of Operations Research, Location Science, and further applied disciplines.
In addition, the most examined concepts of it naturally connect to Physics, Functional Analysis, and non-Euclidean Geometries.
We will not introduce basic notions and terminology of this field going beyond our purpose; for its fundaments the reader is referred to the monograph \cite{tho} and to the surveys \cite{martini 1} and \cite{martini 2}.
The present paper refers to {\em bisectors} in (finite dimensional normed or) \emph{ Minkowski spaces}, i.e., to collections of points which have, in each case, the same distance (with respect to the corresponding norm) to two given points ${\bf x}$, ${\bf y}$ of these spaces. Note that bisectors in Minkowski spaces play an essential role in Discrete and Computational Geometry, mainly in view of constructing (generalized) Voronoi diagrams, and also for motion planning with respect to translations; see, e.g., the surveys \cite{aurenhammer} and \cite{martini 2}.

In some previous papers on this topic (see \cite{gho1}, \cite{gho2}, and \cite{gho3}),
\'A. G. Horv\'ath proved that if the unit ball of a Minkowski space is strictly convex, then every bisector is a topological
hyperplane (see Theorem 2 in \cite{gho1}, or \cite{martini 2}). On the other hand, Example 3 in \cite {gho1} shows that
strict convexity does not follow from the fact that all bisectors are topological hyperplanes.

In these papers, the connections between shadow boundaries of the
unit ball and bisectors (regarding the direction ${\bf x}-{\bf y}$ or the points {\bf x} and {\bf y}, respectively) in Minkowski spaces are investigated. The author was sure
that the following statement is true: \emph{A bisector is a
topological hyperplane if and only if the corresponding shadow
boundary is a topological $(n-2)$-dimensional sphere}. However,
the respective conjecture was proved only in the three-dimensional case (Theorem
2 and Theorem 4 in \cite{gho2}).
In \cite{gho3}, some further topological observations on shadow boundaries are discussed, e.g., that they are
compact metric spaces containing $(n-2)$-dimensional closed,
connected subsets separating the boundary of $K$. \'A. G. Horv\'ath also investigated
the manifold case, and he proved (using an approximation theorem for
cell-like mappings) that the shadow boundary is homeomorphic to an
$(n-2)$-dimensional sphere. Consequently this result (if
the bisector is a homeomorphic copy of $R^{n-1}$, then the shadow
boundary is a topological $(n-2)$-sphere) confirms the
first direction of the above mentioned conjecture.

Independently, H. Martini and S. Wu \cite{martini-wu} introduced and investigated the concept of radial projection of bisectors. Strongly using the central symmetry of Minkowskian balls, they proved some interesting results on radial projections of bisectors.

Theorem 2.6 in \cite{martini-wu} says that the shadow boundary is a subset of the closure of such a radial projection, and
Theorem 2.9 there refers to the converse statement. If for a point ${\bf x}$ from the boundary of the unit ball there exists a point ${\bf z}$, unique except for the sign, such that ${\bf x}$ is orthogonal to ${\bf z}$ in the sense of Birkhoff (see below), then ${\bf z}$ is a point of the radial projection of the bisector corresponding to ${\bf x}$ and $-{\bf x}$.

In the present paper we introduce the concept of bounded representation of bisectors, which yields a useful combination of the notions of bisector, shadow boundary, and radial projection. We prove that the topological properties of the radial projection (in higher dimensions) do not determine the topological properties of the bisector. More precisely, the manifold property of the bisector does not imply the manifold property of the radial projection (see our Example below). The situation is different with respect to the bounded representation of the bisector. Namely, if one of them is a manifold, then the other is also. More precisely, if the bisector is a manifold of dimension $n-1$, then its bounded representation is homeomorphic to a closed $(n-1)$-dimensional ball $B^{n-1}$ (i.e., it is a cell of dimension $n-1$). And conversely, if the bounded representation is a cell, then the closed bisector is also (Theorem 1).

We will also present new approaches to higher dimensional analogues of several theorems given in \cite{martini-wu}. By our new terminology, we will rewrite and reprove Theorems 2.6, 2.9, and 2.10 from that paper.

\section{Basic notions and radial projections of bisectors}

Let $K$ be a compact, convex sets with nonempty interior (i.e., a {\em convex body}) in the $n$-dimensional Euclidean space $E^n$ which, in addition, is centred at the origin $O$. Then the $(n-1)$-dimensional boundary bd $K$ of $K$, in the following also denoted by $S$, can be interpreted as the {\em unit sphere} of an $n$-dimensional (real Banach or) \emph{Minkowski space} $M^n$ with norm $\|\cdot \|$, i.e.,
$$
S:=\{{\bf x}\in M^n\mbox{ : } \|{\bf x}\|=1\}.
$$
 It is well known that there are different types of orthogonality in Minkowski spaces. In particular, for ${\bf x},{\bf y}\in M^n$ we say that ${\bf x}$ is {\em Birkhoff orthogonal} to ${\bf y}$ if $\|{\bf x}+t{\bf y}\|\geq \|{\bf x}\|$ for all $t\in \mathbb{R}$, denoted by ${\bf x}\bot _B{\bf y}$ (see \cite{birkhoff}); and {\bf x} is {\em isosceles orthogonal} to ${\bf y}$ if $\|{\bf x}+{\bf y}\|=\|{\bf x}-{\bf y}\|$, denoted by ${\bf x}\bot _I{\bf y}$ (cf. \cite{james}). The {\em shadow boundary} $S(K,{\bf x})$ of $K$ with respect to the direction {\bf x} is the intersection of $S$ and all supporting lines of $K$ having direction ${\bf x}$.

Given a point ${\bf x}\in
S$, the \emph{bisector} of $-{\bf x}$ and ${\bf x}$, denoted by $B(-{\bf x},{\bf x})$, consists of all those vectors {\bf y}
which are isosceles orthogonal to ${\bf x}$ with respect to the Minkowski norm generated by $K$. The \emph{radial projection} $P({\bf x})$ of this bisector consists of those points ${\bf y}$ of $S$ for which there is a positive real value $t$ such that $t{\bf y}\in B(-{\bf x},{\bf x})$.

We remark that, in the relative topology of $S$, $P({\bf x})$ can either be closed or open; this can be easily seen in the cases of the Euclidean and of the maximum norm. Thus, for
topological investigations in higher dimensions we suggest the extension of the definition of $B(-{\bf x},{\bf x})$ to ideal points by a limit property.

\begin{defi}
Consider the compactification of $E^n$ to a closed ball $B^n$ by the set of ideal points ${\bf x}_\infty $ ($-{\bf x}_\infty\not ={\bf x}_\infty$). We say that ${\bf y}_\infty :=\infty {\bf y}\in B(-{\bf x},{\bf x})$ if  there is a sequence  $(t_i{\bf y}_i)\in B(-{\bf x},{\bf x})$ for which  $\lim\limits_{i\rightarrow \infty} {\bf y}_i={\bf y}$. We call the points of the original bisector {\em ordinary points} and the new points {\em ideal points}, respectively.
\end{defi}

By this definition $P({\bf x})$ could be closed as we can see in our first statement. Let $P({\bf x})^l$ be the collection of those points ${\bf y}$ of $S$ for which $$
\|t{\bf y}+{\bf x}\|< \|t{\bf y}-{\bf x}\|
$$
holds, for all real $t\geq 0$.  Let $P({\bf x})^r$ denote the image of $P({\bf x})^l$ under reflection at the origin.

\begin{proposition}
In the described way, $S$ is decomposed into three disjoint sets: $P({\bf x})$, $P({\bf x})^l$, and $P({\bf x})^r$. $P({\bf x})$ is an
at least $(n-2)$-dimensional closed (and therefore compact) set in $S$ which is connected for $n \geq 3$, the
sets $P({\bf x})^l$ and $P({\bf x})^r$ are arc-wise connected components of their union.
\end{proposition}

\proof By Theorem 5.1 of \cite{martini-wu}, $P({\bf x})$ is connected for $n\geq 3$. We prove that it is also closed with respect to the relative topology of the boundary of the unit ball. To see this, consider a convergent
sequence $({\bf y}_i$) in $P({\bf x})$ having the limit {\bf y}. For any $i$ there is a new sequence of points $({\bf y}^j_i)$ such that for every pair $\{i,j\}$ there are $t_j\in \mathbb{R}^+$ and ${\bf x}_i^j\in  B(-{\bf x};{\bf x})$ such that $ (t^j_i{\bf y}^j
_i) ={\bf x}^j_i$. (For an ordinary point the mentioned sequence can be regarded as a constant one.) It is clear that for the diagonal sequence $({\bf y}^i_i)$ we have
$$
\lim\limits_{i\rightarrow \infty}{\bf y}^i_i = {\bf y},
$$
implying that {\bf y} is also in $P({\bf x})$.
\begin{figure}
    \centerline{\includegraphics[scale=0.8]{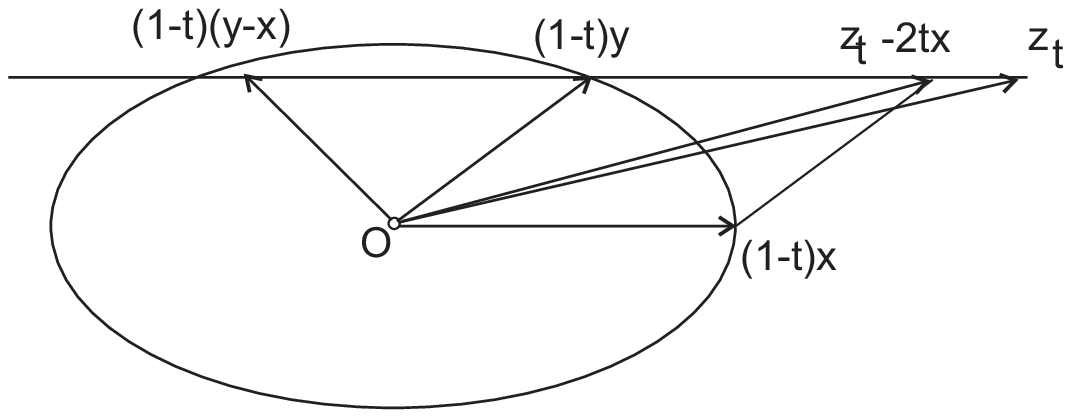}}
\caption{Vectors used in the proof of Proposition 1}
\end{figure}
The continuity property of the norm function implies that all points of $S$ belong to precisely one
of the three mentioned sets. Thus the first statement is clear, and the union of $P({\bf x})^l$ and $P({\bf x})^r$
is open with respect to the topology of $S$.
Observe once more that $P({\bf x})^l$ and $P({\bf x})^r$ are images of each other regarding reflection at the origin. Furthermore, they are arc-wise connected sets. To prove this, consider the following inequality for
an element {\bf y} of $P({\bf x})^r$:
$$
\|({\bf y}-t({\bf y}-{\bf x})-{\bf x}\| = (1-t)\|{\bf y}-{\bf x}\| < (1-t)\|{\bf y} + {\bf x}\| = \|({\bf y}-t({\bf y}-{\bf x})) + {\bf x}- 2t{\bf x}\|,
$$
where $0 \leq t \leq 1$ is an arbitrary parameter. The point ${\bf z}_t := ({\bf y}-t({\bf y}-{\bf x}))+{\bf x} = (1-t){\bf y}+(1+t){\bf x}$
is on the right half-line, starting with the point $(1-t)({\bf y} + {\bf x}) ={\bf z}_t- 2t{\bf x}$ and being parallel to the
vector ${\bf x}$, meaning that its norm is larger than the norm of the point ${\bf z}_t - 2t{\bf x}$ (see Fig. 1).
Thus
$$
\|{\bf z}_t\|\geq \|{\bf z}_t -2t{\bf x}\|,
$$
and so
$$
\|({\bf y}-t({\bf y}-{\bf x}))-{\bf x}\| < \|({\bf y}- t({\bf y}-{\bf x})) + {\bf x}\|.
$$

A consequence of this inequality is that the arc of $S$ connecting the respective endpoints of the
vectors ${\bf y}$ and ${\bf x}$ belongs to the set $P({\bf x})^r$. Thus every two points of $P({\bf x})^r$ can be connected by
an arc, as we stated.
Now, with respect to the topology of their union, they are connected components. This means
that both of them are also open with respect to the topology of $S$. Thus $P({\bf x})$ separates $S$. By
a theorem of P. S. Aleksandrov (Theorem 5.12 in vol. I of \cite{alexandrov}) we get that the topological dimension of
$P({\bf x})$ is at least $n-2$.
\qed

\remark
If we identify the opposite points of $S$, then we get an $(n-1)$-dimensional projective space $P$ dissected into two parts, one of them open and the other one closed, respectively. The open part contains the points of the identified sets $P({\bf x})^r$ and $P({\bf x})^l$, while the closed one contains the identified point-pairs of $P({\bf x})$.

\section{Bounded representation of the bisector}

\begin{figure}
    \centerline{\includegraphics[scale=0.8]{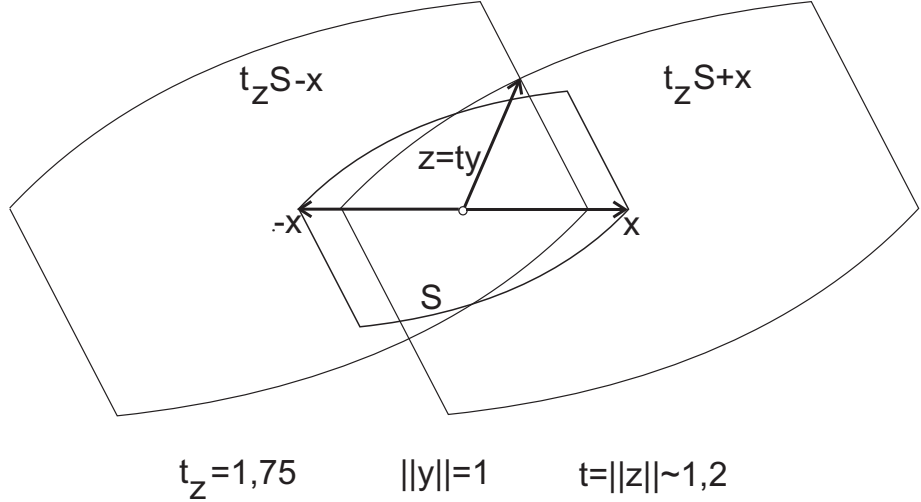}}
\caption{Connection between the parameters $t$ and $t_{\bf z}$.}
\end{figure}

We define now an important mapping of the bisector.
Let ${\bf z}$ be a point of $B(-{\bf x},{\bf x})$. If it is an ordinary point, then there is a unique
value $1 <  t_{\bf z} < \infty$ for which ${\bf z}\in (t_{\bf z}S + x) \cap (t_{\bf z}S-x)$. Let $\Phi  : B(-{\bf x},{\bf x})\longrightarrow K$ denote
the mapping which sends {\bf z} into $\Phi({\bf z}) = \frac{1}{t_{\bf z}}{\bf z}$. For ordinary points the mapping $\Phi $ of the bisector
is continuous. We extend $\Phi $ to the ideal points by the following rule: The image of an ideal point is its radial projection. Denote the
image set of $\Phi $ (with respect to this extended mapping) by $\Phi(B(-{\bf x},{\bf x}))$. We will call this set \emph{the bounded representation of
the bisector}. We now have a connection between the concept of bisector and metrical properties of $K$.

\begin{lemma}
The bounded representation of the bisector is the union of the shadow boundary of $K$ and the locus of the midpoints of the chords of K parallel to {\bf x}.
\end{lemma}

\proof
For an ordinary point {\bf z} of the bisector we have $1 \leq t_{\bf z} < \infty$, and thus the norm of
$$
\frac{1}{t_{\bf z}}{\bf z} =\frac{1}{2}\left( \frac{1}{t_{\bf z}}({\bf z}-{\bf x}) +\frac{1}{t_{\bf z}}
({\bf z} + {\bf x})\right)
$$
is less or equal to 1.

If it is equal to 1, then the point $\frac{1}{t_{\bf z}}{\bf z}$ is a point of a horizontal segment (parallel to {\bf x}) of the boundary and thus a point of the shadow boundary, and the set of all points corresponding to the value $t_{\bf z}$ yields a horizontal segment of $S$. If now $t\geq t_{\bf z}$, the points of the bounded representation corresponding to this value $t$ form
another segment containing the segment of $t_{\bf z}$. Thus the directions determined by the points of the segment of $t_{\bf z}$ are ideal points of the bisector, proving that the points of the shadow boundary are images of certain ideal points.

\begin{figure}
    \centerline{\includegraphics[scale=0.4]{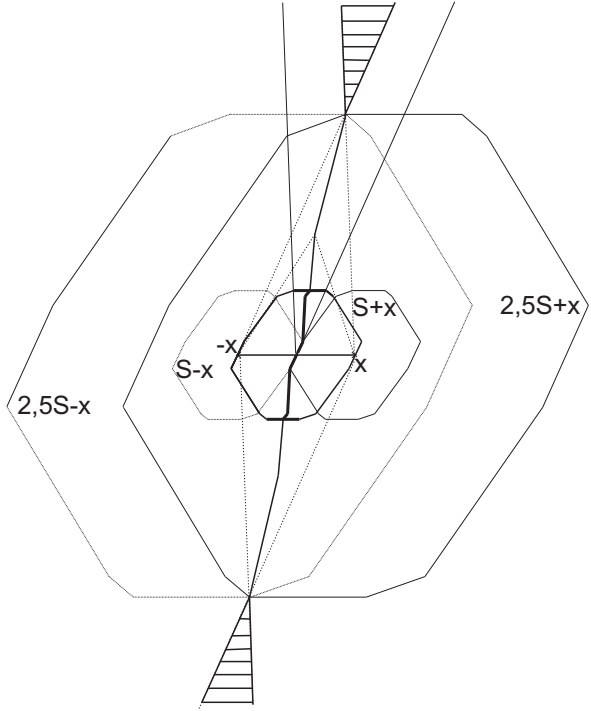}}
\caption{Bounded representation of the bisector}
\end{figure}

In the other case the obtained point is the midpoint of that chord whose endpoints are $\frac{1}{t_{\bf z}}({\bf z}-{\bf x})\in S$ and $\frac{1}{t_{\bf z}}({\bf z}+{\bf x})\in S$, respectively.

Now, by the definition of ideal points, the continuity of the mapping is clear. In fact, we have to check that the image of a point of the bisector with large norm is close to the boundary $S$ of $K$. Since, by definition, $t_{\bf z}$ is equal to $\|{\bf z}-{\bf x}\|$, we have the two inequalities
$$
1\geq \|\frac{1}{t_{\bf z}}{\bf z}\|=\frac{\|{\bf z}\|}{\|{\bf z}-{\bf x}\|}=\frac{1}{\|\frac{{\bf z}}{\|{\bf z}\|}-\frac{{\bf x}}{\|{\bf z}\|}\|}\geq \frac{1}{1+\frac{\|{\bf x}\|}{\|{\bf z}\|}},
$$
showing that for ${\bf z}$ with large norm its bounded representation is close to $S$. To visualize the proof, we show in Fig. 2 the bisector and its bounded representation in a
two-dimensional space.
\qed

For example, by \cite{james2} we get from this lemma immediately the following

\begin{corollary}
{\em The bounded representation of the bisector $B({\bf x},-{\bf x})$ with respect to any point {\bf x} from the unit sphere of a Minkowski space is contained in an $(n-1)$-subspace if and only if the Minkowski space is Euclidean.}
\end{corollary}

It is easy to see that the shadow boundary $S(K,{\bf {\bf x}})$ is the set of points of $S$ which are orthogonal to the vector ${\bf x}$ in the sense of Birkhoff, and thus some results on Birkhoff orthogonality in \cite{martini-wu} can be extended to higher dimensions. These are described in the following.

Theorem 2.6 in \cite{martini-wu} says that the shadow boundary is a subset of the closure of the radial projection. A consequence of the concept of bounded representation of the bisector is the fact that such theorems can be extended to higher dimensions. To prove this, we observe that the radial projection is exactly the radial projection of the bounded representation of the bisector, implied by Lemma 1. So it contains the shadow boundary, extending Theorem 2.6 of \cite{martini-wu}.

Theorem 2.9 in \cite{martini-wu} refers to the converse statement. If for a point ${\bf x}$ of $S$ there exists a unique point ${\bf z}$ (except for
the sign) such that ${\bf x}$ is orthogonal to ${\bf z}$ in the sense of Birkhoff, then ${\bf z}$ is a point of the radial projection of the bisector corresponding to ${\bf x}$. If we denote the sharp points of the shadow boundary as those points of $K$ which are unique in their carrying supporting line of $K$ parallel to ${\bf x}$, we can say that the sharp points of the shadow boundary corresponding to the direction  of ${\bf x}$ belong to the radial projection of the bisector of ${\bf x}$. This is also a consequence of our present definitions and Lemma 1.

Let now $n=2$, and assume that ${\bf z}$ is the radial projection of an ideal point of the bisector of ${\bf x}$. Then it belongs to the shadow boundary of ${\bf x}$ implying that it is Birkhoff orthogonal to ${\bf x}$. Thus, if we consider a boundary point {\bf z} of the radial projection which is not Birkhoff orthogonal to ${\bf x}$, this is a projective image of an ordinary point of the bisector. Thus ${\bf z}$ is a point of the radial projection in the classical sense, too. This proves Theorem 2.10 in \cite{martini-wu}

The main aim of this paper is to prove the following theorem on the topology of the bounded representation of the bisector.

\begin{theorem}
If the bisector is a manifold of dimension $n-1$  with boundary, then its bounded representation is homeomorphic to the $(n-1)$-dimensional closed ball $B^{n-1}$. Conversely, if the bounded representation is a topological ball of dimension $n-1$, then the bisector is of the same type. Furthermore, its relative interior (which is the set of its ordinary points) is a topological hyperplane of dimension $n-1$.
\end{theorem}

\proof Assume that the bisector is a manifold of dimension $n-1$ with boundary. Then an ordinary point has a relatively open $(n-1)$-dimensional neighborhood in the bisector, and thus there are interior points. On the other hand, there is no ideal point which could be in the relative interior of the bisector implying that the set of ordinary points of the bisector is a manifold of dimension $n-1$. Hence our assumption implies that the shadow boundary $S(K,{\bf x})$ is a manifold of dimension $n-2$. In fact, from Theorem 5 and Theorem 4 in \cite{gho3} we get that the shadow boundary is also a topological manifold of dimension $n-2$. Theorem 2 says that it is homeomorphic to $S^{n-2}$. On the other hand, the set $C$ of midpoints of correspondingly directed chords containing interior points of $K$ is always homeomorphic to the positive part $S^+$ of the boundary $S$ of $K$, determined by the shadow boundary. Thus it is homeomorphic to $\mathbb{R}^{n-1}$. Finally we observe that the boundary of the latter set $C$ is the shadow boundary itself, showing that the bounded representation of the bisector is homeomorphic to $B^{n-1}$, as we stated.

We remark that the converse statement is true if and only if the manifold property of the bounded representation can be extended to the bisector. This is clear for the points mapping to the interior of $K$, but it is not evident for other points of the bisector. The problem is that the pre-images of a point of the shadow boundary could form a point or a half-line, respectively. Thus $\Phi $ is not an injective (but, of course, a surjective) continuous mapping. Clearly, both of the two sets (the bisector and its bounded representation) are continua, i.e., compact, connected Hausdorff ($T_2$) spaces. Moreover, the points and half-lines are cell-like sets; thus $\Phi$ is a cell-like mapping. Restricting $\Phi$ to the ideal point of the bisector, we get a bijective mapping onto the shadow boundary. We prove that the set of ideal points is compact in the bisector. Of course, the ideal points of the bisector give a proper part $I$ of $S^{n-1}$ bounding the topological ball $B^n$. Hence this point set can be regarded as a subset of an $(n-1)$-dimensional Euclidean space $\mathbb{R}^{n-1}$. (We can consider ${\bf x}_\infty$ as the center of a stereographic projection.) Its clear that $I$ is bounded. It is also closed by its definition, and so it is compact by the Heine-Borel theorem on compact sets in $\mathbb{R}^{n-1}$ (see, e.g.,  p. 9 in \cite{knopp}). On the other hand, the shadow boundary can also be regarded as an $(n-2)$-sphere embedded into a Euclidean $(n-1)$-space, because ${\bf x}$ is not a point of it. A continuous and bijective mapping from a compact set of $\mathbb{R}^{n-1}$ into $\mathbb{R}^{n-1}$ is a homeomorphism (see again \cite{knopp}). Thus the ideal points of the bisector give a topological $(n-2)$-dimensional sphere.

Now we prove that the ordinary points of the bisector are, with respect to its relative topology, interior points of it. We remark that it is trivial for a point ${\bf z}\in B(-{\bf x},{\bf x})$ if $\Phi ({\bf z})$ is an interior point of $K$, because $\Phi $ (by its definition) is a homeomorphism on the collection of such points  onto the interior of the bounded representation of the bisector. Thus it is also relatively open with respect to the bisector, and this part of the bisector is a topological manifold, homeomorphic to $\mathbb{R}^{n-1}$.

Let now $\Phi ({\bf z})$ belong to the shadow boundary. Since it is a topological sphere of dimension $n-2$, there is a cell of dimension $n-2$ (a homeomorphic copy of a closed ball of dimension $n-2$), namely $Z$,  containing $\Phi ({\bf z})$ in its interior. The pre-image $\Phi^{-1}(\mbox{int }B)$ of the interior $\mbox{int }B$ of $B$ is (by the continuity of $\Phi $) open with respect to the topology of the bisector and contains ${\bf z}$. Thus it has also an interior point with respect to the topology of the bisector.

Finally we observe that from the compactness of $B$ the existence of an $\varepsilon$ follows for which the set
$$
\{{\bf v} \mbox{ : } \|{\bf z}\|-\varepsilon \leq\|{\bf v}\| \leq \|{\bf z}\|-\varepsilon \mbox{, }{\bf v}\in \Phi^{-1}(B)\}
$$
is a closed cone (truncated by two parallel surfaces) containing ${\bf z}$ in its interior. Since the interior of this body is homeomorphic to $\mathbb{R}^{n-1}$, we get that the set of ordinary points is a manifold of dimension $n-1$. In the proof of Theorem 5 from \cite{gho3} it is shown that if the ordinary points of the bisector yield an $(n-1)$-manifold, then it is homeomorphic to $\mathbb{R}^{n-1}$, and Theorem 6 there establishes that it is a topological hyperplane. Thus we proved that the closed bisector is a cell of dimension $n-1$ whose interior can be embedded in the $n$-dimensional Euclidean space in a standard (unknotted) way, as we stated.
\qed

In higher dimensions the topology of the radial projection does (in contrast to the bounded representation) not simplify by the simplification of the topology of the bisector.  Finally we give an example showing that it is possible that the bisector is a manifold with boundary of dimension $n-1$, but $P({\bf x})$ is not a manifold.

\begin{figure}
    \centerline{\includegraphics[scale=0.8]{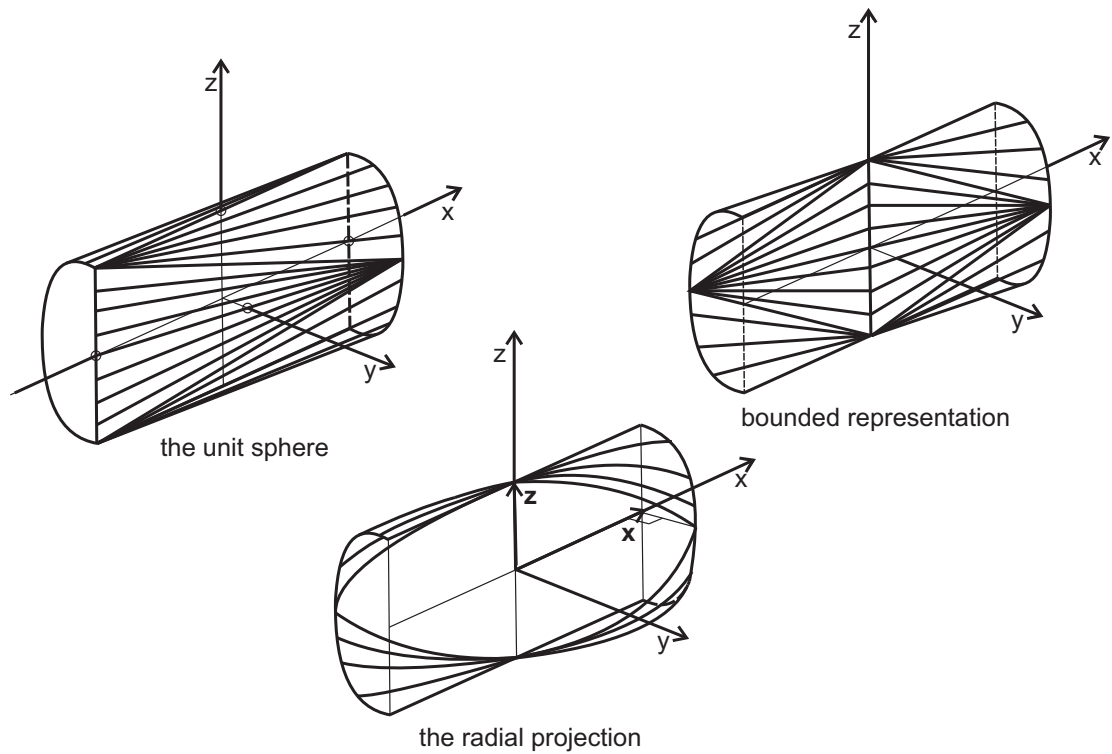}}
\caption{The radial projection is not a manifold}
\end{figure}

\begin{example} Take a cartesian coordinate system in the Euclidean space, with the respective coordinate-axes $x,y$ and $z$, and vectors ${\bf x}$ and ${\bf z}$ having coordinates $(1,0,0)$ and $(0,0,1)$, respectively. Consider the example from \cite{gho1} in which the body $K$ is the convex hull of two half-circles with parallel diameters and in symmetric position with respect to the origin, such that their affine hulls are parallel to the plane $x=0$ (see Fig. 4). Besides the two half circles, the ruled part of the boundary of their convex hull contains four conic surface parts and two opposite triangles. The bounded representation of the bisector corresponding to the direction orthogonal to the planes of the half-circles is homeomorphic to a plane. It can be obtained in the following way: Cut the surface of $K$ by the segments parallel to ${\bf x}$ into two parts. The described half-disks do no longer belong to the surface. Apply an affinity to these two parts, by the ratio $\frac{1}{2}$, with direction orthogonal to the end planes. Finally glue these parts together at their common vertical segment $[-{\bf z},{\bf z}]$. It is clear that the central projection of this ruled surface from the midpoint of $[-{\bf z},{\bf z}]$ is the union of those curves which are the intersections of the planes through $[-{\bf z},{\bf z}]$ and the respective points of $S$. The obtained four parts are conic surfaces separated by $-{\bf z}$, ${\bf z}$ and the two opposite points of the half-circles lying on the horizontal plane $z=0$. Of course, in this case the radial projection is not a topological manifold.
\end{example}

\begin{center}
Horst Martini\\
Faculty of Mathematics\\
Chemnitz University of Technology\\
09107 Chemnitz\\
Germany
\end{center}

\begin{center}
\'Akos G.Horv\'ath,\\
 Department of Geometry \\
Budapest University of Technology and Economics\\
1521 Budapest, Hungary
\\e-mail: ghorvath@math.bme.hu
\end{center}

\end{document}